\documentclass{amsart}
\usepackage[dvips]{graphicx}
\usepackage{amsmath,graphics}
\usepackage{amsfonts,amssymb}

\theoremstyle{definition}

\numberwithin{equation}{section}

\newcommand {\cal } {\mathcal}

\def\Z{{\mathbb Z}}

\def\P{{\mathbb P}}

\def\bn{\begin{enumerate}}
\def\en{\end{enumerate}}
\newtheorem{thr}{Theorem}[section]

\newtheorem{lm}[thr]{Lemma}
\newtheorem{pr}[thr]{Proposition}

\begin{document}

\title[Symplectic representatives of the 
canonical class]
{On the number of components of the symplectic representatives of the 
canonical class}

\author{Stefano Vidussi}
\address{Department of Mathematics, University of California,
Riverside, CA 92521, USA}
\email{svidussi@math.ucr.edu} 
\date{May 13, 2006}
\thanks{The author was supported in part
by NSF grant \#0629956.}


\begin{abstract}
In this paper we show that there exists
a family of simply connected, symplectic $4$-manifolds such that
the (Poincar\'e dual of the)
canonical class admits both connected and disconnected symplectic
representatives. This answers a question raised in \cite{FS4}.
 
\end{abstract}

\maketitle

\section{Introduction and statement of the result}
An important result of $4$-dimensional symplectic topology, proven first by
Taubes in \cite{Ta} using the relation between Seiberg-Witten and Gromov
invariants and then (under minor assumptions) by Donaldson and Smith in \cite{DS}
via Lefschetz fibration techniques, is
the existence of a symplectic representative of the (Poincar\'e dual of the) 
canonical class of a symplectic $4$-manifold with $b_{+}>1$. 
These proofs, in general,
do not provide a sufficiently explicit construction of such a representative,
nor make any statement concerning uniqueness, number of components, 
or their genus. It is therefore a non-obvious task, given a symplectic
$4$-manifold, to provide explicitly such a representative.

An interesting case of this problem is described in \cite{FS4};
Fintushel and Stern show that, for any
choice of positive integers $\{(g_{i},m_{i}), i = 1,...,n\}$,
there exist a (minimal) simply connected symplectic manifold $X$
whose canonical
class $K_{X} \in H_{2}(X,\Z)$ is represented by an embedded symplectic
 surface $\Sigma$ with
$\sum_{i=1}^{n} m_{i}$ connected components:
\begin{equation} \label{rep}
\Sigma = \coprod_{i=1}^{n} \coprod_{j=1}^{m_{i}} \Sigma_{g_{i},j} \in K_{X}, \end{equation}
where $\Sigma_{g_{i},j}$ is a connected surface of genus $g_{i}$.
These manifolds are obtained through natural symplectic operations, i.e.
symplectic fiber sum and symplectic rational blowdown, on simply connected
elliptic surfaces without multiple fibers $E(s)$.
The representative $\Sigma$ of (\ref{rep}) is, in some sense, the
natural result of such operations when we start with the algebraic
representative of the canonical class of the elliptic surface, namely the
disjoint union of $(s-2)$ copies of the fiber.

 Led by this construction,
Fintushel and Stern ask whether, for a symplectic manifold whose canonical class
admits a symplectic representative as in 
(\ref{rep}), the set of integers $\{(g_{i},m_{i}), i = 1,...,n\}$
is a symplectic invariant.
This question is carefully asked under the hypothesis of all $g_{i} \geq 2$;
without this constraint, plenty of counterexamples can be found in \cite{FS3},
\cite{S}, \cite{V}, where it is shown (with different constructions)
that the canonical class of $E(s)$, for $s \geq 4$,
can be symplectically represented by a
connected (non-algebraic) torus.

 In this paper, we will provide an answer, in the negative, to the
aforementioned question, by showing that it is possible to exhibit a
connected symplectic representative for $K_{X}$ for the family of manifolds
constructed in \cite{FS4}. With obvious modifications, symplectic representatives
with any number of components between $1$ and $\sum_{i=1}^{n} m_{i}$ can be obtained.

Roughly speaking, the idea behind our construction consists in ``sewing together" some (or all) the components of the
symplectic representative of (\ref{rep}), while keeping the resulting representative symplectic. For the manifolds $X$ considered in this paper such an internal surgery is explicitly exhibited, and similar cases can be treated analogously. However, it is conceivable that a similar process exists in general, whenever a disconnected representative is available: we are not aware, at this point, of obstructions to the existence of a connected symplectic representative of the canonical class. 

We want to point out that, without contradiction, the result we obtain does not exclude
the use of a numerical symplectic invariant related to the number of components
of symplectic representatives of the canonical class (the set of integers
$\{(g_{i},m_{i}), i = 1,...,n\}$ for a {\it maximal number} of components
$\sum_{i=1}^{n} m_{i}$ could be such an example);
it just stresses the need of
accounting the various representatives. In particular it is possible that the set of integers $\{(g_{i},m_{i}), i = 1,...,n\}$ determined in \cite{FS4} is a symplectic invariant of the family of manifolds therein defined.

{\bf Organization of the paper:} Sections 2 and 3 provide some preliminary material that will be of use in Section 4 for our main construction. More precisely, in Section 2 we discuss a presentation of the elliptic surfaces $E(n)$, 
$n \geq 2$, as symplectic link surgery manifolds, as the first step in identifying some natural submanifolds. In Section 3 we exhibit various symplectic spheres and tori in $E(n)$ that will be the building blocks of our construction. Section 4 contains an inductive presentation of the manifold $X$, reviewing some of the steps of \cite{FS4}, that leads to the construction of the disconnected and connected symplectic representatives of $K_{X}$.

\section{Elliptic surfaces as link surgery manifolds} \label{linksu}
The construction of the manifold $X$ in \cite{FS4} starts by symplectic summing
elliptic surfaces along the fiber $F$ (obviously a symplectic submanifold) and
along a second symplectic torus $R$ (a {\it rim} torus) that arises,
in the surface $E(n+2) = E(n+1) \#_{F_{1} = F_{2}} E(1)$ ($n \geq 0$), by identifying two tori,
in the exterior of the fibers $F_{i}$, that become essential after the sum.
In order to study this construction, we will present an elliptic surface $E(n+2)$
as link surgery manifold (see \cite{FS2} for the definition) obtained from the
Hopf link. This presentation will help us identify some symplectic
submanifolds (spheres and tori) in the elliptic surface $E(n+2)$,
for $n \geq 0$, that we will use in our
construction.

\begin{figure}[h]
\begin{center}
\includegraphics[scale=.2]{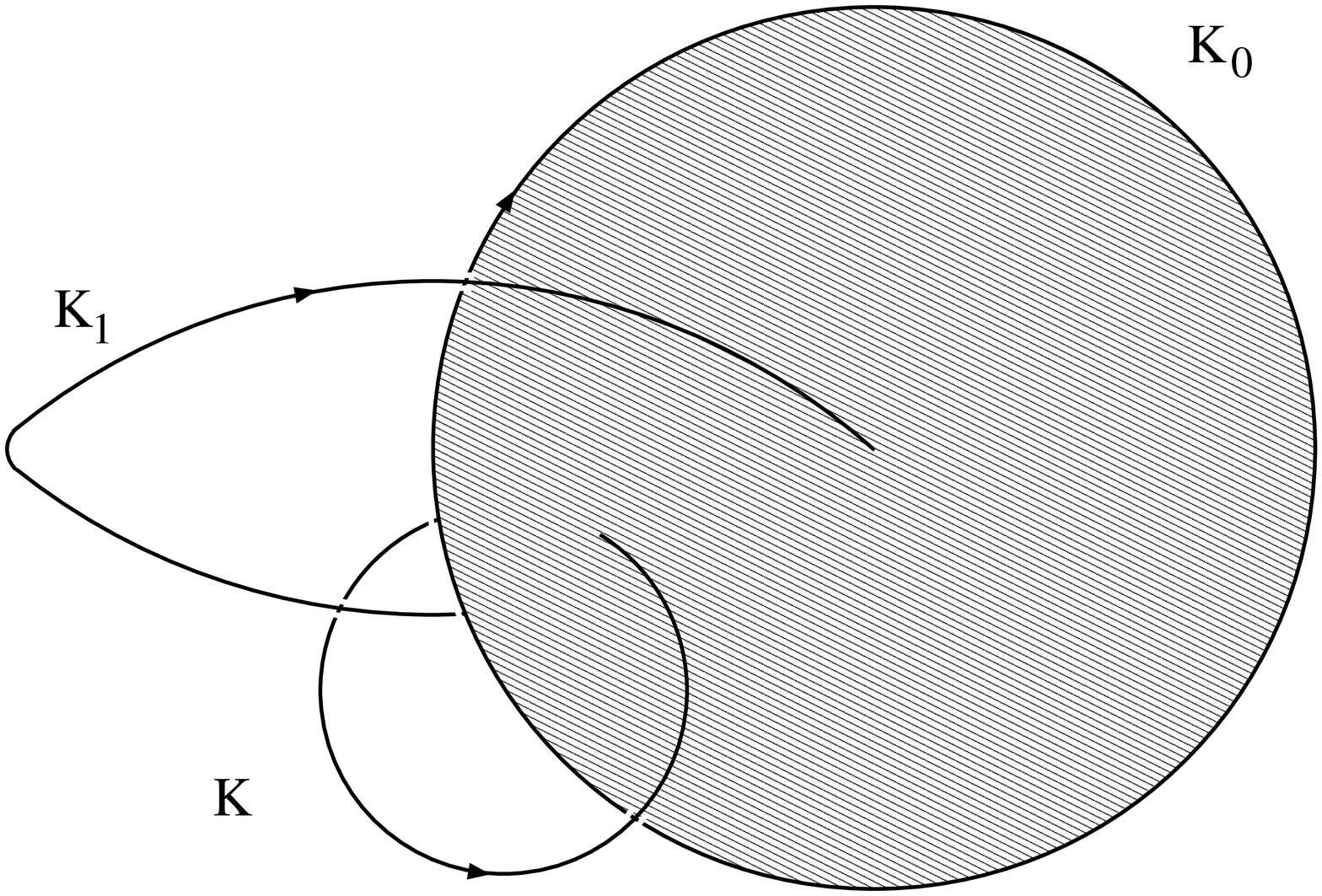}
\end{center}
\caption{The simple closed curve $K$ and the Hopf
link H - the spanning surface of $K_{0}$ is dashed.}
\label{figh}
\end{figure}

Consider the Hopf link $H = K_{0} \cup K_{1}$ and, for future reference,
denote by $K$ the simple closed curve, in $S^{3} \setminus \nu H$, which links
once $K_{0}$ and $K_{1}$ as in Figure \ref{figh}.

\begin{figure}[h]
\begin{center}
\includegraphics[scale=.4]{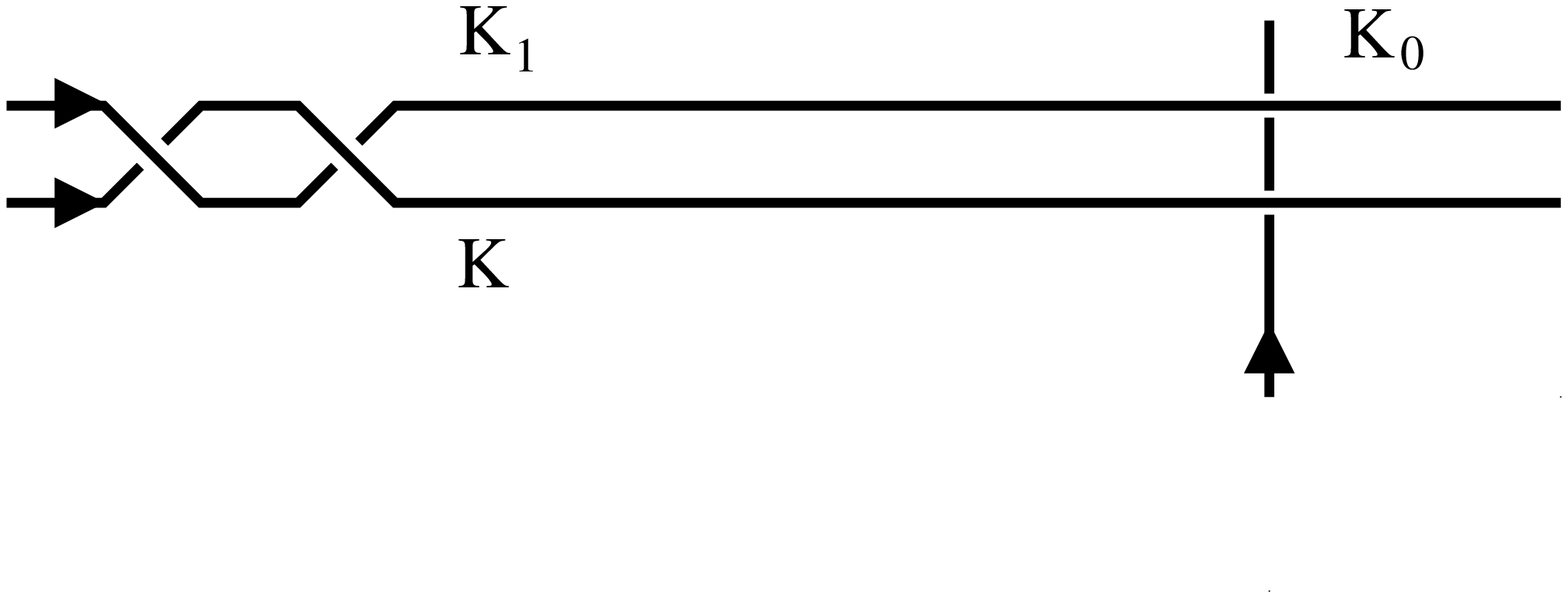}
\end{center}
\caption{The closure of the braid gives $H \cup K$.}
\label{figi}
\end{figure}

A more suggestive presentation of the resulting three component link
$H \cup K$ appears by considering it as the closure, with the axis
originating $K_{0}$, of the braid of Figure \ref{figi}.

The exterior of the link $S^{3} \setminus \nu H$ contains two annuli $A$ and $B$
that are fibers of two distinct fibrations having homology class (in $H^{1}(S^{3} \setminus
\nu H, \Z) = \Z^{2}$) equal to $(1,0)$ and $(0,1)$ respectively. Such annuli have
boundary $\lambda(K_{0}) \cup -\mu(K_{1})$ and $-\mu(K_{0}) \cup \lambda(K_{1})$
respectively. It is useful, for future reference, to think of $B$ as the annulus swept by an arc in the fiber $A$ by the corresponding fibration. 
In what follows, we will consider all the knots endowed with the framing defined
by their spanning disk. Moreover, whenever we will make reference to meridians
and longitudes, we will implicitly assume a specific choice of these
curves is made.

Now observe that any elliptic surface $E(n+2)$
can be presented as the link
surgery manifold obtained by gluing the manifold (with two boundary
components) $S^{1} \times (S^{3} \setminus \nu H)$ to the
exterior of the elliptic fiber of $E(1)$ and $E(n+1)$:
in fact $S^{1} \times (S^{3} \setminus \nu H) = T^{2} \times
(S^{2} \setminus \nu \{p_{0},p_{1}\}) = T^{2} \times A$, where $A$ is
the annulus defined above (removing an open neighborhood of the
Hopf link from $S^{3}$ gives a circle times the annulus $A$,
with the circle identified to the meridian to $K_{0}$).
The usual fiber sum definition of elliptic surfaces can be therefore
interpreted as follows:
\begin{equation} \label{slo} E(n+2) = (E(n+1) \setminus \nu F_{0}) \cup S^{1}
\times (S^{3} \setminus \nu H) \cup (E(1) \setminus \nu F_{1})
\end{equation} where the first gluing map identifies, in the boundary
$3$-tori, $F_{0}$ with $S^{1} \times \mu(K_{0})$ and (remembering that we
reverse orientations) the meridian to $F_{0}$ with $-\lambda(K_{0})$,
while the second gluing map identifies $F_{1}$ with $S^{1} \times \lambda(K_{1})$
and the meridian to $F_{1}$ with $\mu(K_{1})$. After gluing, the fibers $F_{0}$
and $F_{1}$ get identified. Note moreover that
the smooth structure of the resulting manifold is unaffected by the choice
of the diffeomorphism between the fibers of the elliptic surfaces and
$S^{1} \times \mu(K_{0})$ and $S^{1} \times \lambda(K_{1})$; we will later
choose a particular identification.

We need to keep track of the two tori of $E(n+2)$ that
are images of $S^{1} \times \mu(K_0)$ and $S^{1} \times \lambda(K_{0})$ respectively. The first one, identified
with the fiber $F$ of the elliptic fibration of $E(n+2)$, is clearly essential,
but also the second one, that we denote by $R$, is essential in $E(n+2)$,
see e.g. \cite{GS}.
Due to its origin, we will call it rim torus.

To analyze the symplectic submanifolds, we consider how the
construction above leads to present $E(n+2)$ as a symplectic fiber sum:
Perform Dehn surgery along $H \subset S^{3}$ with coefficient $0$ along $K_{0}$
and coefficient $\infty$ along $K_{1}$ to get $S^{1} \times S^{2}$. 
Denote by $C_{0}$ and $C_{1}$ the cores
of the solid tori: these are, up to isotopy, standard circles of the form
$S^{1} \times \{p_{i}\}$ in the
resulting $S^{1} \times S^{2}$. The manifold
$S^{1} \times (S^{1} \times S^{2})$ has a natural symplectic structure
(of the form $dt \wedge \alpha + \epsilon \beta$,
where $\alpha$ represents the fibration of
$S^{1} \times S^{2}$ and $\beta$ is a volume form on the sphere).
The tori $S^{1} \times C_{i}$ are symplectic,
framed, selfintersection zero tori and, after scaling the symplectic forms on each summand if necessary,
we can write \begin{equation} \label{fibsa} E(n+2) =
E(n+1) \#_{F_{0} = S^{1} \times C_{0}} S^{1} \times
(S^{1} \times S^{2}) \#_{S^{1} \times C_{1} = F} E(1).
\end{equation} The symplectic form, away from the gluing locus, restricts
to the symplectic form of the summands.

Note that the presentation of (\ref{fibsa}) can be interpreted, from a
certain viewpoint, as a glorified form of Weinstein's Tubular Neighborhood
Theorem, in the sense that it provides information on the restriction of
the symplectic form of $E(n+2)$ to the submanifold $F \times A$
which separates $E(n+1) \setminus \nu F$ and
$E(1) \setminus \nu F$; application of Weinstein's Theorem to a fiber of $E(1)$,
and then fiber sum with $E(n+1)$, tells that $F \times A$, up to symplectomorphism,
has product symplectic structure (where on the annulus $A$ we take the restriction of a
symplectic form on the sphere). This is the same as the symplectic structure
arising on the image of $S^{1} \times
(S^{1} \times S^{2} \setminus \nu(C_{0} \coprod C_{1})) =
F \times A$ in the fiber sum of (\ref{fibsa}).
\section{Some symplectic spheres and tori in $E(n+2)$} \label{torinel}
In this section we will exploit the presentation of the elliptic surface
$E(n+2)$ of Section \ref{linksu} to identify some symplectically embedded spheres
and tori that will be the building blocks for our (and Fintushel-Stern's) construction. 

We start with the tori. Three symplectic, framed, selfintersection zero tori
arise from the presentations of (\ref{slo}) and (\ref{fibsa}).
The first is the fiber $F$. The second one, $R$, is one
of the two ``marked'' rim tori that become nontrivial in the fiber sum of $E(n+1)$
and $E(1)$ (see Section 3.1 of \cite{GS}),
the other one being identified with the image of $\mu(K_{1}) \times
\lambda(K_{1})$. This essential torus is naturally Lagrangian and becomes symplectic by a small perturbation of the symplectic structure, as discussed in \cite{G}. The third one is the image of $S^{1} \times K$.
Their properties are summarized in the following proposition.
\begin{pr} \label{tori} Consider the inclusion map
\[ S^{1} \times (S^{3} \setminus \nu H)
\hookrightarrow E(n+2). \] The following holds:
\bn 
\item The image $F$ of the torus
$S^{1} \times \mu(K_{0})$ under the inclusion above is a symplectic,
framed, connected submanifold
of $E(n+2)$. 
\item The image $R$ of the torus $S^{1} \times \lambda(K_0)$ is a
Lagrangian, framed, connected submanifold of $E(n+2)$. 
\item The image $T$ of the torus $S^{1} \times K$  is a
symplectic, framed, connected
submanifold of $E(n+2)$ satisfying $[T] = [F] + [R] \in
H_{2}(E(n+2))$. 
\en 
Moreover, the three tori above can be assumed to be 
disjoint.
\end{pr} 
\begin{proof} The first part of the statement
clearly holds true. The non-obvious part is to prove that
$R$ and $T$ are respectively Lagrangian and symplectic w.r.t. the
symplectic structure on $E(n+2)$ induced by the symplectic fiber sum of 
(\ref{fibsa}).
Up to isotopy, we can assume that $\lambda(K_0)$ lies on a spanning disk of 
$K_{0}$ (see Figure \ref{figh}); it is therefore contained in a fiber of the
fibration of $S^{1} \times S^{2}$ obtained by capping off the disk fibration
of $S^{3} \setminus \nu K_{0}$ induced by the spanning disks.
As a consequence, the symplectic structure on $S^{1} \times 
(S^{1} \times S^{2})$  restricts trivially to 
$S^{1} \times \lambda(K_0)$ (its tangent space is spanned by 
$\frac{\partial}{\partial t}$ and a vector tangent to the sphere)
so that the torus $R$ (homologically nontrivial) 
is a Lagrangian submanifold of $E(n+2)$. 
This covers (2) above.
For what concerns (3) we note that, up to isotopy,
$K$ is transversal
to the fibration of $S^{1} \times S^{2}$ (see Figure \ref{figi}),
so that the symplectic form on
$S^{1} \times (S^{1} \times S^{2})$ never vanishes on the torus
$S^{1} \times K$. As this torus is symplectic in one summand,
it will be symplectic in the fiber sum of (\ref{fibsa}).
For what concerns the homology class, observe that, in the homology of $S^{3} \setminus \nu H$,
\begin{equation} \label{homcla}
[K] = \mbox{lk($K$,$K_{0}$)}[\mu(K_{0})] +
\mbox{lk($K$,$K_{1}$)}[\mu(K_{1})] = [\mu(K_{0})] + [\mu(K_{1})]
\in H_{1}(S^{3} \setminus \nu H). \end{equation}
The meridian $\mu(K_{1})$ is homologous (actually, isotopic) to $\lambda(K_{0})$, so the relation follows at this point from the identification of  
$S^{1} \times \mu(K_0)$ and $S^{1} \times \lambda(K_0)$ with $F$ and $R$ respectively.
The fact that these tori are disjoint follows directly from the
construction. \end{proof}

The output of the previous proposition, namely that the (primitive) class
$[F] + [R]$ can be represented by two symplectic submanifolds, one given by
the disjoint union $F \coprod R$ and the second by $T$ will be, in fact,
the main tool for our construction. Note that the curve $K$ can be interpreted as result of \textit{circle summing} the meridian $\mu(K_{0})$ and the longitude $\lambda(K_{0})$. The resulting operation on the tori $F$ and $R$, that produces the symplectic torus $T$ from the symplectic torus $F$ and the Lagrangian torus $R$, represents at local level the ``sewing" referred to in the Introduction.

\begin{figure}[h]
\begin{center}
\includegraphics[scale=.3]{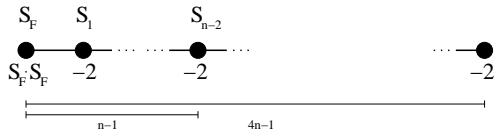}
\end{center}
\caption{The first $n-1$ spheres give the configuration
$\Gamma_{n}$ in the elliptic surface $E(n+2)$.} 
\label{diag}
\end{figure}

Having dealt with tori, we will now consider spheres. Specifically, we are
interested in two groups of spheres. The first group are ``sections" of elliptic nuclei, where an \textit{elliptic nucleus} is the regular neighborhood of the union of a cusp fiber and a section of an elliptic fibration. The second group of spheres are those contained in a configuration that we can use for \textit{rational blowdown}. Remember that this surgery consists in replacing, in a $4$-manifold, a regular neighborhood of a configuration $\Gamma_n$ of $n-1$ spheres as in Figure \ref{diag}, the first with self--intersection $-(n+2)$ and the remaining of self--intersection $-2$, with a certain rational homology ball $B_n$. This rational homology ball naturally embeds in the Hirzebruch surface ${\bf F}_{g-1}$, as the
exterior of the configuration of spheres $(S_{+} +f) \cup S_{-}$, where $S_{+}$ ($S_{-}$) is the positive (negative) 
section and $f$ the fiber of the sphere fibration of ${\bf F}_{g-1}$; see \cite{FS1} for the details of this construction.   

We will start with the first group. It is well known that an elliptic surface
$E(n+2)$
contains several disjoint elliptic nuclei (see \cite{GS} for example).
We will be interested in two of them. We have the following proposition:
\begin{pr} Let $E(n+2)$, $n \geq 0$ be an elliptic surface, with the symplectic structure inherited from the construction of Section \ref{linksu}. The following holds true: 
\bn 
\item $E(n+2)$ contains an elliptic nucleus $N_{F}$, with symplectic fiber $F$ and symplectic sphere $S_{F}$ of selfintersection $-(n+2)$, given by a section of the elliptic fibration. 
\item $E(n+2)$ contains an elliptic nucleus $N_{R}$, with Lagrangian fiber $R$ and Lagrangian sphere $S_{R}$ of selfintersection $-2$; with a small perturbation of the symplectic structure of $E(n+2)$ we can make $R$ and $S_{R}$ symplectic.
\item The nuclei $N_{F}$ and $N_{R}$ are disjoint,
and the torus $T$ intersects their spheres in a single, transverse point.
\en \end{pr}
\begin{proof}  The first nucleus $N_{F}$ arises, in the picture of 
(\ref{fibsa}), as regular neighborhood
of a cusp fiber of the elliptic fibration and its $-(n+2)$-sphere
section $S_{F}$; the symplectic fiber of this nucleus
is $F$, and the section $S_{F}$ arises by gluing together a disk
section of
$E(n+1) \setminus \nu F$ with a $-(n+1)$ twist {\it rel $\partial$} in its normal
bundle (for sake of brevity we will denote these disks, with the usual abuse of notation,
$-(n+1)$-disks), a $(-1)$-disk section of $E(1) \setminus \nu F$ and
(for a suitable choice of the embedding in
$S^{1} \times (S^{3} \setminus \nu H) = F \times A$)
the annulus $A \subset S^{3} \setminus \nu H$ (which has boundary identified with
$\lambda(K_{0}) \cup  - \mu(K_{1})$). The two vanishing disks that kill
the generators of $\pi_{1}(F)$ can be located in $E(n+1) \setminus \nu F_{0}$, as
$F$ is already contained in a nucleus in $E(n+1)$. $S_{F}$ is symplectic, as it
appears as connected sum of symplectic spheres in each summand of 
(\ref{fibsa}). This proves (1).

The second nucleus $N_{R}$ contains, as regular fiber, the
Lagrangian rim torus $R$. The vanishing disks that kill the generators of
$\pi_{1}(R)$ can be located in
$E(1) \setminus \nu F_{1}$ (using a second disk section and a vanishing disk).
The $(-2)$-sphere of the nucleus $S_{R}$ is obtained,
in the construction of (\ref{slo}), by gluing the annulus $B$
(which has boundary identified with $-\mu(K_{0}) \cup \lambda(K_{1})$) embedded in
$S^{1} \times (S^{3} \setminus \nu H) =
F \times A$ to two $(-1)$-vanishing disks in
$(E(n+1) \setminus \nu F_{0})$ and $(E(1) \setminus \nu F_{1})$.
The annulus $B$ is Lagrangian, as we can span at each point its tangent space
with a vector $v$ in $ker \alpha$ and a vector $w$ that satisfies
$\beta(w,\cdot)=0$, so that the symplectic form on $F \times A$ vanishes
on $B$. Similarly, the vanishing disks in the elliptic surfaces
can be taken to be Lagrangian thimbles of a symplectic Lefschetz fibration (see \cite{ADK}, Section $4$ and, for a general discussion,
\cite{D}):
we start with a generic pencil of cubics in $\P^{2}$, $\pi: \P^{2} \setminus B
\rightarrow \P^{1}$ (where the base locus $B$ is composed of $9$ points) and
we endow it with a symplectic connection, given by the symplectic orthogonal
to the fiber.
Then, given a path $\gamma(t)$ on $\P^{1}$ joining a critical value $\pi(p_{c}) =
\gamma(0)$ with a
regular value $\gamma(1)$ (and otherwise disjoint from the set of critical
values) we can define a Lagrangian vanishing disk as the union
of the vanishing cycles on the fibers lying over the
path, defined by the condition that symplectic parallel transport ${\cal P}$
sends them to
the critical point on the singular fiber, namely
\begin{equation} p_{c} \cup \bigcup_{t} V_{t} := p_{c} \cup \bigcup_{t} \{ u
\in \pi^{-1}(\gamma(t)) |\lim_{\epsilon \rightarrow 0}
{\cal P}_{\gamma |_{[\epsilon,t]}} u = p_{c} \}.
\end{equation}
When we symplectically blowup $\P^{2}$ along the base locus of the pencil to
get $E(1)$, we obtain from the exceptional divisors $9$ symplectic sections,
disjoint from the Lagrangian disk. Further fiber summing, to get other elliptic
surfaces, does not affect the vanishing disk.
When we recover $E(n+2)$ through the fiber sum of (\ref{fibsa})
we choose the identification map on the two boundary tori in such a way
to identify the boundary of the annulus $B$ with the boundaries of the
vanishing disks, obtaining this way the Lagrangian $(-2)$-sphere $S_{R}$.
The essential Lagrangian submanifolds $R$ and
$S_{R}$ can be made symplectic with a small perturbation of the symplectic form, as in
\cite{FS4}. This covers (2).

For what concerns (3), note that the
intersection of $A$ and $B$
in $S^{3} \setminus \nu H$ gets removed in
$F \times A$, due to the presence of the extra $S^{1}$ factor.
Moreover, thanks to the abundance of singular fibers and sections,
we can choose all the disks used in the previous construction to be
disjoint (see \cite{GM}). As a result, the nuclei $N_{F}$ and $N_{R}$ are
disjoint.
Finally, the torus $T$ intersects $S_{F}$ and $S_{R}$ only in the interior
of $F \times A$, and the intersection is a single transverse point corresponding
to the intersection of $K$ with the annuli $A$ and $B$, as it appears
from Figure \ref{figh}. \end{proof}

\begin{figure}[h]
\begin{center}
\includegraphics[scale=.5]{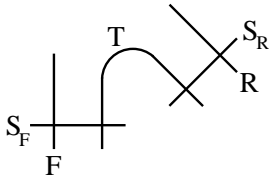}
\end{center}
\caption{Schematic representation of the relation between the two nuclei of
fiber $F$ and $R$ and the torus $T$. }
\label{scheme}
\end{figure}

Figure \ref{scheme} schematizes the relation between the two nuclei and the
torus $T$.

The second group of spheres we are interested in is the
configuration of symplectically embedded spheres $\Gamma_n$.
As discussed in \cite{FS1} any elliptic surface
$E(n)$ contains a configuration of $4n - 1$
symplectic spheres as in the diagram of Figure \ref{diag}.

The sphere of
selfintersection $-n$ is the section of the elliptic fibration,
and the configuration intersects the fiber $F$ only in $S_{F}$,
with a single transverse intersection point.
When we fiber sum $E(2)$ to $E(n)$ (for $n \geq 1$)
we have therefore (Lemma 2.1 of \cite{FS4}), in the resulting $E(n+2)$,
a configuration
of $4n-1$ embedded symplectic spheres as in the diagram of Figure \ref{diag}
(where the sphere $S_{F}$ of selfintersection $-(n+2)$, obtained
by gluing two sections, is again a section of the elliptic fibration).
In particular, if we keep track only of the first $(n-1)$ spheres of the
configuration, we deduce that $E(n+2)$, for $n \geq 1$, contains the configuration $\Gamma_{n}$, 
whose spheres will be denoted as
$S_{F},S_{1},...,S_{n-2}$.
The same statement holds true if we
consider the manifold (with boundary) $N(2) \#_{F} E(n)$.
Note that in $E(n+2) = E(2) \#_{F} E(n)$,
as the nucleus determined by the rim torus $R$
is disjoint from the section of $E(2)$, it is disjoint also from
the configuration $\Gamma_{n}$; moreover the torus $T$
(entirely contained in $E(2) \setminus \nu F$) intersects $\Gamma_{n}$
in a single transverse point on the first sphere $S_{F}$.

\section{Fintushel-Stern construction and connected representatives}
For sake of clarity, we summarize some of the results discussed in
the previous two sections. In the elliptic surface $E(n+2)$,
for $n \geq 1$, we can identify two disjoint symplectically
embedded
surfaces. The first is the linear plumbing of the $(n-1)$ spheres
$S_{F},S_{1},...,S_{n-2}$ described in the diagram of Figure \ref{diag}.
The second is the nucleus $N_{R}$. The fiber $F$ intersects
$S_{F}$ in a single positive point, orthogonal w.r.t. the symplectic
structure, while it does not intersect the other
spheres of $\Gamma_{n}$ nor the nucleus $N_{R}$. The symplectic torus $T$,
instead, intersect both $S_{F}$ and $S_{R}$ in
a single transverse point, and is otherwise disjoint from the remaining
spheres and tori considered above.

In the situation above (except for the torus $T$, that has no role
in their construction) Fintushel and Stern have
inductively defined a family of minimal symplectic simply connected 
$4$-manifolds whose canonical class admits a symplectic representative 
satisfying the condition of (\ref{rep}). 
We will sketch their construction and
show how, by a suitable modification, we can obtain a connected representative 
for the canonical class. (In what follows we will reserve the symbol ${\hat \cdot}$ to denote connected
surfaces.)

First, for $g \geq 2$, let $X_{g}$ be the
manifold obtained by symplectic rational blowdown of the configuration
$\Gamma_{g}$ in $E(g + 2)$, i.e. replacing $\nu \Gamma_g$ with the rational homology ball $B_g \subset {\bf F}_{g -1}$ (see Section \ref{torinel}). The resulting manifold is simply connected and symplectic, where
the symplectic structure is obtained by grafting to $E(g+2) \setminus \nu
\Gamma_{g}$ the symplectic structure induced by the embedding of
$B_{g}$ in the Hirzebruch surface ${\bf F}_{g -1}$ endowed with a 
suitable symplectic structure (see \cite{Sy}). 
A symplectic representative of the canonical class $K_{X_{g}}$ 
(image of $K_{E(g+2)} = g[F]$ under the blowdown map) is then given by gluing 
$(F_{1} \coprod ... \coprod F_{g}) \setminus \nu \Gamma_{g}$ (a collection of $g$
copies of the fiber with a hole) to $S_{+} \cap B_{g}$ (a sphere with 
$g$ holes, as $S_{+}$ and $(S_{+} +f) \cup S_{-}$ intersect in $g$ points).
The result of this surgery is therefore a connected, embedded, symplectic 
surface of genus $g$ that we will denote by ${\hat \Sigma_{g}}$, 
and which represents
$K_{X_{g}}$. 

There is a homology class we want to keep track 
of, and which represents the image, under the blowdown, of the class
$g[F] + [R]$. The nucleus $N_{R}$ is disjoint from
$\Gamma_{g}$, so that it survives the blowdown process; we index with an $R_X$
its image.
(A more correct, but notationally heavier, index would be $R_{X_{g}}$.)
The image of the class $g[F] + [R]$ is therefore well
defined and is given by $[{\hat \Sigma_{g}}] + [R_{X}] = K_{X_{g}} +
[R_{X}]$. The disjoint union
${\hat \Sigma_{g}} \coprod R_{X}$ is a symplectic representative, with two
connected components of genus $g$ and $1$, but (as $[T] = [F] + [R]$)
we can represent is as well
by the connected, symplectic, genus $g$ surface
${\hat \Sigma}_{g}^{R_{X}}$ obtained by gluing, much as above,
$(F_{1} \coprod ... \coprod F_{g-1} \coprod T)
\setminus \nu \Gamma_{g}$
(a collection of $g$ disjoint tori with a hole) to
$S_{+} \cap B_{g}$. The surface ${\hat \Sigma}_{g}^{R_{X}}$
intersects the $(-2)$-sphere $S_{R_{X}}$ in a single positive transverse
point, which is the image of the intersection point of $T$ and $S_{R_{X}}$
under the blowdown.

Summing up, we have a symplectic manifold $X_{g}$ which has
a canonical class $K_{X_{g}}$ represented by ${\hat \Sigma_{g}}$,
a connected symplectic surface of genus $g$. This manifold contains an
embedded symplectic nucleus $N_{R_{X}}$ (with fiber $R_{X}$) disjoint from
${\hat \Sigma_{g}}$; moreover the class $K_{X_{g}} + [R_{X}]$ can be
represented
by the disjoint union ${\hat \Sigma_{g}} \coprod R_{X}$ or by a connected 
symplectic surface ${\hat \Sigma}^{R_{X}}_{g}$ of genus $g$.

The manifold $X_{g}$ is the initial step in the inductive construction
of the family of manifolds of \cite{FS4}, and now we will proceed to the
inductive step. In practical terms, we need to go over the proof of Lemma 2.2 of \cite{FS4} and ensure that there is room to use, at each inductive step, a torus of the type $T$ to ``sew" the components of the disconnected representative identified by Fintushel and Stern. The argument is quite straightforward but rather long to present.

\begin{lm} \label{conc} Let $\{g_{i},i=1,...,m\}$ a collection of integers $g_{i} \geq 2$. Let $X$ be a symplectic simply connected $4$-manifold satisfying the following conditions:
\bn
\item 
The canonical class $K_{X}$ of $X$ can be represented by the
union $\Sigma_{g_{1},...,g_{m}}$ of
$m$ disjoint connected symplectic surface of genus $g_{1},...,g_{m}$ or by a
connected symplectic surface ${\hat \Sigma}_{g_{1},...,g_{m}}$ of genus $(\sum_{1}^{m}g_{i} - m +1)$.
\item $X$ contains a symplectic
nucleus $N_{R_{X}}$ with fiber $R_{X}$ and section $S_{R_{X}}$ disjoint from
both $\Sigma_{g_{1},...,g_{m}}$ and ${\hat \Sigma}_{g_{1},...,g_{m}}$.
The class $K_{X} + [R_{X}]$ can be represented by the union
of $m+1$ disjoint connected symplectic surfaces $\Sigma_{g_{1},...,g_{m}} \coprod R_{X}$
or by a connected symplectic
surface ${\hat \Sigma}^{R_{X}}_{g_{1},...,g_{m}}$ of
genus $(\sum_{1}^{m}g_{i} - m +1)$ intersecting the $(-2)$-sphere
$S_{R_{X}}$ in a single positive transverse point. 
\en
Then for any $g \geq 2$, there is a symplectic simply connected
manifold $Y$ satisfying the following conditions:  
\bn 
\item The canonical class $K_{Y}$ of $Y$ can be represented by the
union $\Sigma_{g_{1},...,g_{m},g}$ of
$m+1$ disjoint connected symplectic surface of genus $g_{1},...,g_{m},g$ or
by a connected symplectic surface ${\hat \Sigma}_{g_{1},...,g_{m},g}$ of genus
$(\sum_{1}^{m}g_{i} + g - m)$. 
\item Y contains a symplectic
nucleus $N_{R_{Y}}$ with fiber $R_{Y}$ and section $S_{R_{Y}}$ disjoint from both
$\Sigma_{g_{1},...,g_{m},g}$ and ${\hat \Sigma}_{g_{1},...,g_{m},g}$.
The class $K_{Y} + [R_{Y}]$ can be represented by the union
of $m+2$ disjoint connected symplectic surfaces $\Sigma_{g_{1},...,g_{m},g} \coprod R_{Y}$ or by a connected symplectic surface
${\hat \Sigma}^{R_{Y}}_{g_{1},...,g_{m},g}$
of genus $(\sum_{1}^{m}g_{i} +g - m)$ intersecting the $(-2)$-sphere
$S_{R_{Y}}$ in a single positive transverse point. 
\en
\end{lm}

\begin{proof} Following \cite{FS4}, we observe that along the symplectic torus
$R_{X} \subset X$ we can define the symplectic
fiber sum \begin{equation} \label{iden} X \#_{R_{X} = F} E(g) \end{equation}
where $F$ is the standard fiber in $E(g)$. The resulting manifold
is simply connected and symplectic, with canonical class $K_{X\#_{R_X=F}E(g)} = K_{X} + K_{E(g)} + 2[F] = K_{X} + g[F]$. We have two symplectic representatives for
such class. The first is the disjoint union
$\Sigma_{g_{1},...,g_{m}} \coprod F_{1} \coprod ... \coprod F_{g}$, where the $F_{j}$'s are
parallel copies of $F$; this surface has $m + g$ connected
components. For the second one, observing that
the sum of (\ref{iden}) identifies $R_{X}$ and $F$, we can
choose the disjoint union ${\hat \Sigma}^{R_{X}}_{g_{1},...,g_{m}} \coprod F_{2} \coprod ...
\coprod F_{g}$; this surface has $g$ components.
Figure \ref{case} schematizes the situation for $X \#_{R_{X} = F} E(4)$.

\begin{figure}[h]
\begin{center}
\includegraphics[scale=.5]{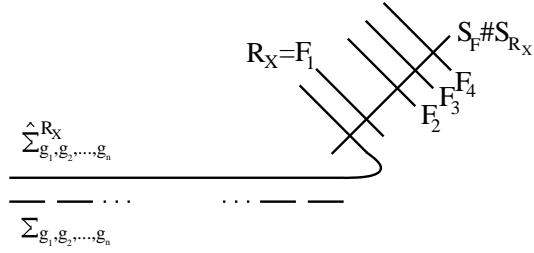}
\end{center}
\caption{Schematic presentation of the components of the two representatives
of $K_{X} + g[F]$ (for $g = 4$); the first one is
$\Sigma_{g_{1},...,g_{m}} \coprod_{i=1}^{4} F_{i}$ and the second one is
${\hat \Sigma}^{R_{X}}_{g_{1},...,g_{m}} \coprod_{i=2}^{4} F_{i}$.}
\label{case}
\end{figure}

As $g \geq 2$, the elliptic surface $E(g)$ of fiber $F$
contains a rim torus $R$ and a torus $T$ satisfying the conditions of
Proposition \ref{tori}. The homology class $K_{X\#_{R_X=F}E(g)} + [R]
= K_{X} + g[F] + [R]$ is well defined and can be
symplectically represented by the disjoint union $\Sigma_{g_{1},...,g_{m}}
\coprod F_{1} \coprod ... \coprod F_{g} \coprod R$,
composed by $m+g+1$ components, as
well as by the disjoint union ${\hat \Sigma}^{R_{X}}_{g_{1},...,g_{m}}
\coprod F_{2}
\coprod ... \coprod F_{g-1} \coprod T$, composed of $g$ components.
As in homology we have $[T] = [F] + [R]$,
these two surfaces are homologous.

Next, we observe that the
manifold $X \#_{R_{X} = F} E(g) \supset N_{R_{X}} \#_{R_{X} = F} E(g)$ contains a
symplectically embedded configuration $\Gamma_{g}$, inherited from a configuration in $E(g)$ that intersects $F$ in a single point (its $-(g+2)$-sphere is the
connected sum of $S_{R_{X}}$ and $S_{F}$).
This configuration is disjoint from
$\Sigma_{g_{1},...,g_{m}}$
(as the $(-2)$ sphere $S_{R_{X}}$ is disjoint from it) and from the symplectic
nucleus with fiber $R$. Instead it does intersect the surface
${\hat \Sigma}^{R_{X}}_{g_{1},...,g_{m}}$ in a single positive transverse
point in the $-(g+2)$-sphere.
We blow down the configuration $\Gamma_{g}$. Denote by $Y$ the
resulting simply connected, symplectic manifold. We point out that the nucleus
$N_{R}$ survives the blowdown.

 The canonical class of $Y$
is given by the image of $K_{X \#_{R_X=F} E(g)}$ under the blowdown, and we will use
the previous constructions for exhibiting two symplectic representatives, one
having $m+1$ connected components of genus $g_{1},...,g_{m},g$, the other one
connected and of genus $(\sum_{1}^{m}g_{i} + g - m)$.

First the disconnected representative $\Sigma_{g_{1},...,g_{m},g}$
(presented in \cite{FS4}) is
obtained by the disjoin union of $\Sigma_{g_{1},...,g_{m}}$ (unaffected by the
blowdown) and a genus $g$ connected surface ${\hat \Sigma_{g}}$ obtained,
as in
the initial step of our construction, by gluing the $g$ tori with hole
$(F_{1} \coprod ...\coprod F_{g}) \setminus \nu \Gamma_{g}$ to the sphere with $g$ holes
$S_{+} \cap B_{g}$. The connected representative
${\hat \Sigma}_{g_{1},...,g_{m},g}$, instead, is obtained
by gluing to $S_{+} \cap B_{g}$ the surface
$({\hat \Sigma}^{R_{X}}_{g_{1},...,g_{m}} \coprod F_{2} \coprod ... \coprod
F_{g}) \setminus \nu \Gamma_{g}$. The resulting surface
is clearly homologous to the previous one, and is connected, as
both $F$ and ${\hat \Sigma}^{R_{X}}_{g_{1},...,g_{m}}$ intersect,
in a positive transverse point,
the $-(g+2)$-sphere of $\Gamma_{g}$. The genus of this surface is then easily
computed. In reference to the scheme of Figure \ref{case}, the two
representatives above
are obtained blowing down $\Gamma_{g}$, of which $S_{R_{X}} \# S_{F}$ is the first
sphere, and connecting the surfaces hit by that sphere. This completes
the proof of (1). 

Next, we observe that the manifold $Y$ contains a symplectic nucleus,
inherited from the one of $E(g)$, whose image we index by $R_Y$. This
nucleus is disjoint, by construction, from both
$\Sigma_{g_{1},...,g_{m},g}$ and ${\hat \Sigma}_{g_{1},...,g_{m},g}$.
To finish our argument, we must proceed to identify a disconnected and
a connected symplectic representative of the class $K_{Y} + [R_{Y}]$,
as stated in (2). For what concerns the
disconnected representative, this is simply provided by the
disjoint union of the disconnected symplectic surface
$\Sigma_{g_{1},...,g_{m},g}$ and a copy of $R_{Y}$. In order to obtain
the connected representative, we consider, in $X \#_{R_X=F} E(g)$, the
surface ${\hat \Sigma}^{R_{X}}_{g_{1},...,g_{m}} \coprod F_{2}
\coprod ... \coprod F_{g-1} \coprod T$, composed of $g+1$ components. By
blowing down, as ${\hat \Sigma}^{R_{X}}_{g_{1},...,g_{m}}$, $F_{i}$ and
$T$ intersect the configuration $\Gamma_{g}$ only in
one point (on the sphere $S_{F}$), the class $K_{Y} + [R_{Y}]$
has a symplectic connected representative
${\hat \Sigma}^{R_{Y}}_{g_{1},...,g_{m},g}$ obtained by gluing,
much as above, the genus $(\sum_{1}^{m}g_{i} - m +1)$ surface
${\hat \Sigma}^{R_{X}}_{g_{1},...,g_{m}}$ with
one hole, and the $(g-1)$ tori with hole
$F_{2} \coprod ... \coprod F_{g-1} \coprod T$ to the sphere with $g$ holes
$S_{+} \cap B_{g}$. A check, using the genus formula for connected sum, or the
adjunction formula, shows that the genus of the surface is the one stated.
This surface intersects $S_{R_{Y}}$ in a single positive transverse
point, which is the image of the intersection point of $T$ and $S_{R}$
under the blowdown. \end{proof}

Lemma \ref{conc} provides the inductive step required, and completes the construction of a symplectic connected
surface, homologous to the disconnected symplectic representative of (\ref{rep}).

\end{document}